\documentclass[11pt]{article}
\usepackage{amsthm} 
\usepackage{mathrsfs}
\usepackage{amsmath}
\usepackage{amsfonts}
\usepackage{amssymb, amsmath, cite}
\usepackage{color}
\usepackage{graphicx}

\newtheorem{theorem}{Theorem}[section]
\newtheorem{problem}{Problem}[section]

\renewcommand{\theequation}{\thesection.\arabic{equation}}

\newcommand\be{\begin{equation}}
\newcommand\ee{\end{equation}}
\newcommand\ber{\begin{eqnarray}}
\newcommand\eer{\end{eqnarray}}
\newcommand\berr{\begin{eqnarray*}}
\newcommand\eerr{\end{eqnarray*}}
\newcommand\lm{\lambda}
\newcommand\lb{\label}\newcommand\eq{\eqref}

\newcommand{\dd}{\mathrm{d}}\newcommand\bfR{\mathbb{R}}\newcommand\bfF{\mathbb{F}}\newcommand\bfC{\mathbb{C}}
\newcommand\pa{\partial}
{\setlength\arraycolsep{2pt}
\newcommand\nn{\nonumber}
\newcommand\ch{\mathrm{ch}}
\newcommand\bea{\begin{eqnarray}}
\newcommand\eea{\end{eqnarray}}

\usepackage{graphicx}
\setlength{\baselineskip}{18pt}{\setlength\arraycolsep{2pt}

\begin{document}

\title{The Jacobian Conjecture and Integrability of\\ Associated Partial Differential Equations\footnote{Invited contribution for a memorial issue of Physica D in honor of Vladimir Zakharov}}

\author{Yisong Yang\,\footnote{Email address: yisongyang@nyu.edu}\,\,\footnote{ORCID number: 0000-0001-7985--6298}\\Courant Institute of Mathematical Sciences\\New York University}

\date{}

\maketitle

\begin{abstract}

The Jacobian conjecture over a field of characteristic zero is considered directly in view of the nonlinear partial differential equations it 
is associated with. 
Exploring the integrals of such partial differential equations, this work obtains broad families of polynomial
maps satisfying the conjecture in all dimensions and of arbitrarily high degrees. Furthermore, it is shown that a reformulated multiply parametrized
version of the conjecture in all dimensions enables a separation of the Jacobian equation into a system of
subequations which may be integrated systematically rendering a settlement of the parametrized Jacobian problem in this context.

\medskip

\medskip


\end{abstract}

\section{Introduction}
\setcounter{equation}{0}

In this article, we consider the Jacobian conjecture and relate its resolution, either way, to the polynomial integrability of the families of
fully nonlinear partial differential equations it is associated with.

Let $\bfF$ be a field of characteristic zero, $\ch(\bfF)=0$, and consider a polynomial map
\be
P:\bfF^n\to\bfF^n,\quad n\geq1.
\ee
That is $P=P(x)=(P_1(x),\dots,P_n(x))$ where $P_1,\dots,P_n$ are some polynomials in $x=(x_1,\dots,x_n)\in\bfF^n$. We are interested in the condition under which $P$ has a polynomial
inverse. That is, there is a polynomial map $Q:\bfF^n\to\bfF^n$, $Q=Q(y),y\in\bfF^n$, such that $P(Q(y))=y$ and $Q(P(x))=x$ for all $x,y\in\bfF^n$. From this, if $\bfF=\bfR$ or $\bfC$, we can compute the Jacobian matrices of $P$ and $Q$, say $\dd P (x)$ and $\dd Q(y)$, to get 
\be
\dd P(x)\dd Q(y)=\dd Q(y)\dd P(x)=I,
\ee
which leads to the necessary condition 
\be\lb{03}
J(P)(x)\equiv \det(\dd P(x))=a\in\bfF,\quad a\neq0, 
\ee
for the Jacobian of the map $P$, since both $J(P)(x)$ and $J(Q)(y)$ are polynomials in $x$ and $y$, respectively. When $n=1$, the condition \eq{03} indicates that $P$ is simply of the
form $P(x)=ax+b$ ($a\neq0$) such that $Q(y)=\frac{y-b}a$, where $a,b\in\bfF$ and $\bfF$ can be any field with zero characteristic. It  is easy to see that the condition \eq{03} is not sufficient to ensure
the invertibility of $P$ when $\ch(\bfF)\neq0$. For example, consider the polynomial $P(x)=x-x^p$ for $x\in \bfF$ and $p=\ch(\bfF)\geq2$. Then $\dd P(x)=1$ but $P$ fails to be 
one-to-one even.

The Jacobian conjecture was first formulated in 1939 by Ott-Heinrich Keller \cite{K} which states that a polynomial
map $P:\bfF^n\to \bfF^n$ where $\bfF$ is a field of characteristic zero and $n\geq2$ has a polynomial inverse if the Jacobian
$J(P)$ of $P$ is a nonzero constant. This conjecture has not been solved for any
$n\geq2$ and appears in Smale's list of
the eighteen mathematical problems for the new century \cite{S}. With a suitable normalization, that is, an affine linear transformation of variables, the polynomial
map $P$ may be taken to satisfy the condition
$
P(0)=0$ and $DP(0)=I$ such that it has the representation
\be\lb{1.1}
P(x_1,\dots,x_n)=(x_1+H_1,\dots,x_n+H_n),\quad (x_1,\dots,x_n)\in\bfF^n,
\ee
where $H_1,\dots,H_n$ are polynomials in the variables $x_1,\dots,x_n$ consisting of terms of degrees at least 2 in nontrivial or nonlinear situations so that
the condition imposed on $J(P)$ becomes $J(P)=1$.
Among the notable developments, Wang \cite{W} established the conjecture when $H_1,\dots,H_n$ are all quadratic  and Bass, Connell, and Wright \cite{BCW} and Yagzhev \cite{Y} proved an important reduction theorem which states that the general conjecture
amounts to showing that the conjecture is true for the special case when each of $H_1,\dots,H_n$ is either cubic-homogeneous or zero
for {\em all} $n$. Subsequently, Druzkowski \cite{D} further showed that the 
cubic-homogeneous reduction of \cite{BCW,Y} may be
assumed to be of the form of cubic-linear type,
\be\lb{1.2}
H_i=(a_{i1}x_1+\cdots+a_{in}x_n)^3,\quad i=1,\dots,n.
\ee
In \cite{BE}, Bondt and Essen proved that the conjecture for the case $\bfF=\bfC$ may be reduced to showing that the conjecture is true when the Jacobian matrix of the map $H=(H_1,\dots,H_n)$ is homogeneous, nilpotent, and symmetric, for all $n\geq2$.
For $n=2$, Moh \cite{M} established the conjecture when the degrees of $H_1$ and $H_2$ are up to 100.
See the survey articles \cite{Dru,E1,Mei} and monograph \cite{E2} and references therein for further results and progress. While these developments
were mainly based on ideas and methods of algebra and algebraic geometry,  the problem also naturally prompts us to
explore its structure in view of partial differential equations directly. With this approach, in the first part of this contribution,
we will present some new families of polynomial maps satisfying the conjecture, which may be summarized as follows.

\begin{theorem}\label{th1.1}
In its simplest form, the Jacobian conjecture amounts to solving an under-determined first-order nonlinear partial differential equation and obtaining relevant polynomial-function solutions which
may be used to construct the inverse map from the original polynomial map over $\bfF^n$ for a field $\bfF$ of characteristic zero.
\begin{enumerate}

\item[(i)] When $n=2$, the equation may be reduced into a homogeneous Monge--Amp\`{e}re equation over $\bfF^2$ whose homogeneous solutions give rise to a broad family of solutions
of any prescribed degree and the inverse polynomial map may be readily constructed based on an invariance structure of the variables so that
the independent and dependent variables of the map are related by the same linear combinations of the variables
and the degrees of the map and its inverse coincide.

\item[(ii)] When $n=3$, the equation may be reduced into an under-determined equation of the Monge--Amp\`{e}re type given in terms of some Hessian determinants involving three unknown
functions. For this equation, broad families of polynomial solutions of arbitrary degrees can be constructed among which, one family is of homogeneous type involving two arbitrary polynomial functions
and a single invariance structure so that the degree of the map and its inverse coincide as in the $n=2$ situation, and another is not of the homogeneous type but involves two arbitrary polynomial functions in composition and obeys a partial invariance structure so that the degree of the inverse map may be as high as twice of that of the original map.

\item[(iii)] For any $n\geq2$,  by imposing a full invariance structure for the variables, a broad family of solutions involving $n-1$ arbitrary polynomial functions
of any prescribed degrees of a single linear combination of
the variables may be constructed explicitly so that the degrees of the map and its inverse coincide.
\end{enumerate}
\end{theorem}

We note that, if $P:\bfF^n\to \bfF^n$ ($n\geq2$) is a polynomial automorphism (that is, the inverse $P^{-1}$ of $P$ exists
which is also a polynomial map), then it has been shown \cite{BCW,Dru,RW} that there holds the following general bound between the degrees of $P$
and $P^{-1}$:
\be\lb{2.13}
\deg(P^{-1})\leq (\deg(P))^{n-1}.
\ee
In the $n=2$ situation here, the map $P$ stated in (i) of Theorem \ref{th1.1} satisfies $\deg(P)=\deg(P^{-1})$ where $\deg(P)$
may be any positive integer.  Besides, the results in (ii) of Theorem \ref{th1.1} says $P^{-1}$ may
fulfill both $\deg(P^{-1})=\deg(P)$ and $\deg(P^{-1})=(\deg(P))^2$ where $\deg(P)$ can be any positive integer as well. In other words, the bound \eq{2.13} is seen to be sharp
or realizable when $n=2,3$ in all degree number cases, with the solutions obtained here.

From the viewpoint of partial differential equations, the main difficulty with the Jacobian conjecture lies in the fact that the normalized Jacobian equation, $J(P)=1$, is a single equation,
which is vastly underdetermined with mixed nonlinearities involving the derivatives of the unknowns. On the other hand, the reduction, or more precisely, separation, method effectively employed in the first part of this study impels us to find
a natural reduction mechanism for the method. For this purpose,
in the second part of this study, we modify the polynomial map \eq{1.1} into
\be\lb{1.4}
P(x_1,\dots,x_n)=(\lm_1 x_1+H_1,\dots,\lm_n x_n+H_n),\quad (x_1,\dots,x_n)\in\bfF^n,
\ee
where $\lm_1,\dots,\lm_n$ are taken to be {\em free} parameters
for the purpose of an effective splitting of the single Jacobian equation. It is clear that if \eq{1.4} has a polynomial inverse for fixed $\lm_1,\dots,\lm_n$, then
\be\lb{1.5}
J(P)=\lm_1\cdots\lm_n,
\ee
and $\lm_1,\dots,\lm_n$ satisfy
\be\lb{1.6}
\lm_1\neq0,\quad \dots,\quad \lm_n\neq0.
\ee
Thus, conversely, we ask whether any polynomial solution to \eq{1.5} has a polynomial inverse under the condition \eq{1.6}. We refer to this statement  as the {\em parametrized Jacobian conjecture} or problem and
we settle this conjecture affirmatively and constructively. More precisely, we will prove: 

\begin{theorem}\lb{th1.2} For the multiply parametrized polynomial map \eq{1.4} where $H_1,\dots,H_n$ are polynomials comprised of quadratically and more highly powered  terms of the variables, all the solutions to the Jacobian equation \eq{1.5}  can be schematically and explicitly constructed. Besides,
under the condition \eq{1.6}, these obtained solutions render the map \eq{1.4} invertible in the category of polynomial maps and the associated inverse maps may also be
constructed schematically and explicitly.
\end{theorem}

In Sections \ref{sec2}--\ref{sec5},  we establish Theorem \ref{th1.1}. In Section \ref{sec6}, we establish Theorem \ref{th1.2}. Since being able to construct explicitly the
associated inverse
polynomial maps is of obvious importance in applications, in Section \ref{sec7}, we briefly illustrate an application to cryptography. 
In Section \ref{s8}, we present the multiparametrized Jacobian problem in its general setting and draw conclusion in such a context.
In Section \ref{sec8}, we summarize the results of this work.

\section{Two dimensions}\lb{sec2}
\setcounter{equation}{0}

First we consider $n=2$ and rewrite \eq{1.1} conveniently as 
\be\lb{2.1}
P(x,y)=(x+f(x,y),y+g(x,y)),
\ee
where $f$ and $g$ are polynomials in the variables $x,y$ over $\bfF$ consisting of terms of degrees at least 2
in nontrivial situations. Inserting \eq{2.1}
into $J(P)=1$ we have
\be\lb{2.2}
f_x+g_y+J(f,g)(x,y)=0,
\ee
where $f_x$ (e.g.) denotes the partial derivative of $f$ with respect to $x$ and $J(f,g)(x,y)$
the Jacobian of the map $(f,g)$ over $x,y$. That is, $J(f,g)(x,y)=\frac{\pa(f,g)}{\pa(x,y)}$. This is an underdetermined equation which may be solved by
setting
\be\lb{2.3}
f_x+g_y=0,\quad J(f,g)(x,y)=0,
\ee
separately, such that the first equation in \eq{2.3} implies that there is a polynomial $h(x,y)$ serving as a scalar potential of the divergence-free vector field $(f,g)$
satisfying
\be\lb{2.4}
f=h_y,\quad g=-h_x.
\ee
Inserting \eq{2.4} into the second equation in \eq{2.3} we see that $h$ satisfies the homogeneous Monge--Amp\`{e}re equation
\cite{A,E}
\be\lb{2.5}
\det(D^2 h)=h_{xx}h_{yy}-h^2_{xy}=0.
\ee
Alternatively, if we are only concerned with $f$ and $g$ being homogeneous of the same degree, then
a degree counting argument applied to \eq{2.2} leads to
two separate equations, as given in \eq{2.3}, as well. Hence we arrive at \eq{2.5} again.

It is of interest and relevance to briefly elaborate a bit more on the Monge--Amp\`{e}re equation in the context of its broader background.

{\bf Gauss curvature and Monge--Amp\`{e}re equation}

Consider a surface in $\bfR^3$ represented as the graph of a function $z=u(x,y)$. Then the Gauss curvature of the surface is given by
\be\lb{02.6}
K=\frac{u_{xx}u_{yy}-u_{xy}^2}{(1+u_x^2+u_y^2)^2}.
\ee
An interesting question in geometric analysis asks whether any prescribed function $K=K(x,y)$ may be realized as the Gauss curvature of a certain surface,
a version of the Nirenberg problem \cite{An,A,Yang}. In view of \eq{02.6}, it is
seen that the question amounts to solving the differential equation of the type:
\be\lb{02.7}
\det(D^2 u)=F(x,y,u,Du),
\ee
classically known \cite{CH} as the Monge--Amp\`{e}re equation over $\bfR^2$. The equation is classified into elliptic, hyperbolic, and homogeneous types, according to $F>0, F<0$, and $F=0$, respectively.
The Jacobian problem, \eq{2.5},  that is, $F\equiv0$ in \eq{02.7}, corresponds to the zero Gauss curvature situation in  \eq{02.6}.  Even in this simple situation, a full integration of
the equation is not yet available and known solutions of the equation assume rather complicated forms. For example, if we take the ansatz
\be\lb{02.8}
u(x,y)=(a_1 x+b_1y+c_1)\sigma\left(\frac{a_2 x+b_2y+c_2}{a_3 x+b_3y+c_3}\right)+a_4 x+b_4y+c_4,
\ee
where $a,b,c$'s are arbitrary constants and $\sigma(\xi)$ is an arbitrary function of $\xi$, then we have
\be
\det(D^2u)=-\frac{(\sigma'(\xi)\det(A))^2}{(a_3 x+b_3 y+c_3)^4},\quad \xi=\frac{a_2 x+b_2y+c_2}{a_3 x+b_3y+c_3},\quad A=\left(\begin{array}{ccc}a_1&b_1&c_1\\a_2&b_2&c_2\\a_3&b_3&c_3\end{array}\right).
\ee
Hence, \eq{02.8} is a solution to the homogeneous Monge--Amp\`{e}re equation if and only if $\det(A)=0$. Geometrically, this implies that a surface of zero Gauss curvature may look surprisingly 
complicated. In \cite{BGZ}, the question of integrability of the equation
\be
\det(D^2u)=c,
\ee
where $c$ is a constant, is investigated through constructing the Lax pair representations. See also \cite{Ban,Gut,KK,Kh,Pav}.
For our interest here, that is, the polynomial solutions, we still do not know what the most general solutions to \eq{2.5} should look like in such a context.

We further note that the Monge--Amp\`{e}re equation and its generalized versions, in $n$ dimensions ($n\geq2$) and with extended nonlinearities, also arise in vast areas of applications including
various geometric analysis problems \cite{Cheng,TW},  optimal mass transport problems \cite{Caf,Cullen}, geometric structures and interpretation of fluid dynamics equations \cite{Banos,Na,Rou}, optics \cite{Brix,Mer},  string and brane theories \cite{Bah,FY,HT,S1,S2} in high-energy physics, and extended gravity models \cite{Ry,Sol}.

We now return to \eq{2.5}. Note that, for our interest here, we are only interested in polynomial solutions with nonlinear terms of sufficiently high powers to be made specific soon. Thus, from \eq{02.8}, we see that the solutions of \eq{2.5} 
that we are looking for are necessarily of the homogeneous type:
\be\lb{2.6}
h(x,y)=\sigma(\xi),\quad \xi=ax+by,\quad a,b\in\bfF,
\ee
as suggested by \eq{1.1}--\eq{1.2}, since we need to take $a_3=b_3=0$ and $c_3\neq0$ which imposes the condition
\be
\left|\begin{array}{cc} a_1&b_1\\a_2&b_2\end{array}\right|=0,
\ee
so that the factor in front of $\sigma$ in \eq{02.8} may be absorbed into the argument of $\sigma$, where the arbitrary polynomial function $\sigma(\xi)$ is taken to be of degree $m\geq3$ or zero. Thus, with the
notation $P(x,y)=(u,v)$ and the relations \eq{2.1} and \eq{2.4}, we have
\be\lb{2.7}
u=x+b\sigma'(\xi),\quad v=y-a\sigma'(\xi),
\ee
resulting in the invariance condition between the two sets of the variables:
\be\lb{2.8}
au+bv=ax+by=\xi.
\ee
This key structure, just unveiled by the Monge--Amp\`{e}re equation \eq{2.5},  will be exploited effectively in all higher dimensions as well.
As a consequence of \eq{2.7} and \eq{2.8}, we obtain the inverse of the map $P$ immediately as follows:
\be\lb{2.9}
x=u-b\sigma'(\xi),\quad y=v+a\sigma'(\xi),\quad\xi=au+bv.
\ee
As a by-product, the arbitrariness of the function $\sigma$ indicates that the solution gives rise to a family of
polynomial maps of arbitrarily high degrees.

For later development, we also observe that the second equation in \eq{2.3} implies that $f$ and $g$ are
functionally dependent. Therefore, if we set $g=G(f)$ (say), then the first 
equation in \eq{2.3} leads to
\be
f_x+G'(f)f_y=0,
\ee
which has a nontrivial solution of the homogeneous type, $f=\phi(\xi)$ ($\xi=ax+by$), if and only if
\be\lb{2.11}
g=G(f)=-\frac ab f,\quad b\neq0.
\ee
 Thus, it follows that there hold the simplified relations
\be\lb{2.12}
u=x+\phi(\xi),\quad  v=y-\frac ab \phi(\xi),\quad au+bv=ax+by=\xi,
\ee
where the invariance relation between the variables again makes the inverse of the map ready to be read off.

We remark that \eq{2.12} is the most general polynomial automorphism of homogeneous type in two dimensions. To see this, we let the map $P$ be defined by
\be\lb{2.14}
u=x+\phi(\xi),\quad v=y+\psi(\eta),\quad \xi=ax+by,\quad \eta=cx+dy,\quad a,b,c,d\in\bfF,
\ee
 where $\phi(\xi)$ and $\psi(\eta)$ are polynomials in the variables $\xi$ and $\eta$, respectively, of
degrees $l\geq2$ and $m\geq2$, satisfying $\phi(0)=\psi(0)=0$. Inserting \eq{2.14} into $J(P)=1$, or $f=\phi$ and
$g=\psi$ into \eq{2.2}, we get
\be\lb{2.15}
a\phi'(\xi)+d\psi'(\eta)+(ad-bc)\phi'(\xi)\psi'(\eta)=0.
\ee
On the other hand, as polynomials in the variables $x,y$, we have
\bea
\deg(a\phi'(\xi)+d\psi'(\eta))&\leq& \max\{l-1,m-1\},\lb{2.16}\\
 \deg(\phi'(\xi)\psi'(\eta))&=&(l-1)+(m-1)>\max\{l-1,m-1\},\lb{2.17}
\eea
since $l,m\geq2$. In view of \eq{2.15}--\eq{2.17}, we arrive at $ad-bc=0$. In other words, the variables
$\xi$ and $\eta$ as given in \eq{2.14} are linearly dependent. Consequently, \eq{2.14} is simplified into
the form
\be\lb{2.18}
u=x+\phi(\xi),\quad v=y+\psi(\xi),\quad \xi=ax+by,
\ee
which renders $a\phi'(\xi)+b\psi'(\xi)=0$. Thus  \eq{2.12} follows if $b\neq0$ and $\phi(\xi)=-\frac ba\psi(\xi)$ if $a\neq0$. So we have obtained the most general homogeneous solution to the equation
$J(P)=1$ or \eq{2.2}.

It can be checked directly that, when the polynomial map \eq{2.1} is of the type $\deg(P)\leq3$, then 
$J(P)=1$ or \eq{2.2} leads to the homogeneous form \eq{2.18}, or more precisely, \eq{2.12}.

One may wonder whether \eq{2.3} is too strong a condition for the splitting of the single equation \eq{2.2}. Here we remark that it arises naturally in three dimensions from a reduction
consideration. In fact, in this situation, we may consider the polynomial map
\be
P(x,y,z)=(x+\zeta(z)f(x,y),y+\zeta(z)g(x,y),z),
\ee
where $\zeta(z)$ is an arbitrary polynomial of $z$. Then  $J(P)=1$ gives us
\be\lb{2.19}
\zeta(z)(f_x+g_y)+\zeta^2(z) (f_x g_y-f_y g_x)=0,
\ee
leading to \eq{2.3} again. Therefore, with the solution \eq{2.6} and notation $P(x,y,z)=(u,v,w)$, the expression \eq{2.7} is updated with
\be
u=x+b\zeta(z)\sigma'(\xi),\quad v=y-a\zeta(z)\sigma'(\xi),\quad w=z,
\ee
so that the invariance property \eq{2.8} still holds. As a consequence, we obtain the inverse map
\be
x=u-b\zeta(w)\sigma'(\xi),\quad y=v+a\zeta(w)\sigma'(\xi),\quad \xi=au+bv,\quad z=w,
\ee
immediately, which extends the formula \eq{2.9}.

This study prompts us to consider the following integrability problems in the context of the Jacobian conjecture.

\begin{problem}\lb{P2.1}
 Can the homogeneous Monge--Amp\`{e}re equation \eq{2.5} be integrated completely in the context of its polynomial solutions? That is, can one obtain all polynomial
solutions of the equation?
\end{problem}

Note that, if $\bfF=\bfR$, the solutions to \eq{2.5} are classified as representing various types of developable surfaces. Thus, for our interest of polynomial
solutions, we are inevitably led to \eq{2.6}. On the other hand, however, if $\bfF$ is an arbitrary field of characteristic zero, this problem remains open.

\medskip

More generally, we may also consider the Jacobian equation \eq{2.2} directly as follows.

\begin{problem}\lb{P2.2}
 Can the Jacobian equation \eq{2.2} be integrated completely in the context of its polynomial solutions? In other words, can one obtain all polynomial
solutions of this equation?
\end{problem}

The study of these problems is important for the resolution of the Jacobian conjecture in two dimensions, in either way. For example, to construct
 a counterexample to the conjecture, one needs to obtain a peculiar solution to \eq{2.5} or \eq{2.2} such that the associated map \eq{2.1} does not
have a
polynomial inverse.

The sense of the integrability of these problems is narrower than that of the usual integrability problem of a system of differential equations in the sense of Liouville (cf. Arnold \cite{Ar}).
Note that, in the study of differential equations, sometimes it is useful to pursue integrability in a specific sense (e.g., the Puiseux integrability \cite{Dem}
and the Weierstrass integrability \cite{Gine}). For our interest in the Jacobian problem, it is only of
relevance to study the polynomial integrability of the differential equations as stated. Of course, our construction yields us abundant families of non-polynomial solutions as well, giving rise to
invertible maps satisfying the Jacobian equation too.

It will also be of interest to know whether Problem \ref{P2.1} and Problem \ref{P2.2} are actually related in the context of polynomial solutions. Specifically,
we ask:

\begin{problem}\lb{P2.3} In the context of polynomial solutions, does the equation \eq{2.2} always split into the two equations in \eq{2.3}?
\end{problem}

An affirmative answer to this problem implies that Problems \ref{P2.1} and \ref{P2.2} are equivalent.

\section{Three dimensions}\lb{sec3}
\setcounter{equation}{0}

Next we consider $n=3$ and rewrite \eq{1.1} as
\be\lb{3.1}
(u,v,w)=P(x,y,z)=(x+f(x,y,z),y+g(x,y,z),z+h(x,y,z)),\quad(x,y,z)\in\bfF^3,
\ee
where $f,g$, and $h$ are polynomials in $x,y,z$ with terms of degrees at least 2
in nontrivial situations. Thus the equation
$J(P)=1$ is recast into
\be\lb{3.2}
f_x+g_y+h_z+J(f,g)(x,y)+J(g,h)(y,z)+J(f,h)(x,z)+J(f,g,h)(x,y,z)=0,
\ee
which is under-determined as well. If we focus on $f,g$, and $h$ being homogeneous of the same degree, then
\eq{3.2} splits into the coupled system
\be\lb{3.3}
f_x+g_y+h_z=0,\quad J(f,g)(x,y)+J(g,h)(y,z)+J(f,h)(x,z)=0,\quad J(f,g,h)(x,y,z)=0,
\ee
as in Section 2. However, here, we are interested in solutions of more general characteristics.

To proceed, we solve the third equation in \eq{3.3} by setting $h=H(f,g)$ where $H$ is a function of 
the variables $f$ and $g$ to be determined. Hinted by the study in Section 2, we seek for solutions of the form
\be\lb{3.4}
f(x,y,z)=\phi(\xi),\, g(x,y,z)=\psi(\eta),\, \xi=ax+by+cz,\, \eta=px+qy+rz,\, a,b,c,p,q,r\in\bfF.
\ee
Inserting \eq{3.4} into the first equation in \eq{3.3}, we get
\be\lb{3.5}
(a+cH_f)\phi'(\xi)+(q+rH_g)\psi'(\eta)=0.
\ee
We are interested in being able to allow $\phi$ and $\psi$ to be arbitrary. This leads to $a+cH_f=0$
and $q+rH_g=0$ or
\be\lb{3.6}
h=H(f,g)=-\frac ac f-\frac qr g,\quad c,r\neq0,
\ee
which extends \eq{2.11}.
In view of \eq{3.6}, we see that the second equation in \eq{3.3} is equivalent to the equation
\be\lb{3.7}
(ar-cp)(br-cq)=0.
\ee
Thus either $ar=cp$ or $br=cq$. 
In other words, subject to \eq{3.7}, we have solved the Jacobian equation
$J(P)=1$ where $P(x,y,z)=(u,v,w)$ in 3 dimensions with
\be\lb{3.9}
u=x+\phi(\xi),\,v=y+\psi(\eta),\, w=z-\frac ac\phi(\xi)-\frac qr \psi(\eta),\,\xi=ax+by+cz,\,\eta=px+qy+rz.
\ee

First, suppose
\be\lb{3.8a}
\frac ac=\frac pr,\quad \frac bc\neq \frac qr.
\ee
Using \eq{3.8a}, we see that \eq{3.9} gives us
\bea
au+bv+cw&=&ax+by+cz+\left(b-\frac{cq}r\right)\psi(\eta)=\xi+\left(b-\frac{cq}r\right)\psi(\eta),\lb{3.10}\\
 pu+qv+rw&=&px+qy+rz=\eta.\lb{3.11}
\eea
So the quantity $\eta$ is seen as an invariant between the two sets of the variables but not $\xi$. In
other words, we achieve a {\em partial invariance}.
As a consequence of \eq{3.9}--\eq{3.11}, we obtain the inverse of the map $P$ given by
\bea
x&=&u-\phi\left(au+bv+cw-\left[b-\frac{cq}r\right]\psi(\eta)\right),\lb{3.12}\\
y&=&v-\psi(\eta),\lb{3.13}\\
 z&=&w+\frac ac\phi\left(au+bv+cw-\left[b-\frac{cq}r\right]\psi(\eta)\right)+\frac qr \psi(\eta),\lb{3.14}
\eea
where $\eta=pu+qv+rw$. 

Next, similarly assume
\be\lb{3.8b}
\frac ac\neq \frac pr,\quad \frac bc= \frac qr.
\ee
Then \eq{3.9} and \eq{3.8b} lead to
\bea
au+bv+cw&=&ax+by+cz=\xi,\\
pu+qv+rw&=&px+qy+rz+\left(p-\frac{ar}c\right)\phi(\xi)=\eta +\left(p-\frac{ar}c\right)\phi(\xi).
\eea
Thus $\xi$ is an invariant but not $\eta$. This again realizes a partial invariance and gives rise to the inverse map analogously by the expressions
\bea
x&=&u-\phi(\xi),\lb{3.18b}\\
y&=&v-\psi\left(pu+qv+rw-\left[p-\frac{ar}c\right]\phi(\xi)\right),\lb{3.19b}\\
z&=&w+\frac ac\phi(\xi)+\frac qr\psi\left(pu+qv+rw-\left[p-\frac{ar}c\right]\phi(\xi)\right),\lb{3.20b}
\eea
where $\xi=au+bv+cw$.

It will be of interest to compare the degrees of the map $P$ given by \eq{3.9} and its inverse $P^{-1}$
either given by \eq{3.12}--\eq{3.14} or \eq{3.18b}--\eq{3.20b}, with regard to the general bound \eq{2.13}, which are
\be
\deg(P)=\max\{\deg(\phi),\deg(\psi)\};\quad \deg(P^{-1})=\deg(\phi)\deg(\psi).
\ee
Hence we have 
\be\lb{deg}
\deg(P^{-1})\leq (\deg(P))^2,
\ee
 which is a realization of \eq{2.13} when $n=3$. Of course,
$\deg(P^{-1})=(\deg(P))^2$ if and only if $\deg(\phi)=\deg(\psi)$ and a wide range of integer combinations 
in the inequality \eq{deg} can
be achieved concretely by choosing appropriate pair of the generating polynomials, $\phi$ and $\psi$.

It is worth noting that, if both factors in \eq{3.7} vanish, or
\be\lb{3.15}
\frac ac= \frac pr,\quad \frac bc=\frac qr
\ee
are simultaneously valid, then \eq{3.10} and \eq{3.11} imply that both $\xi$ and $\eta$ are invariant quantities between
the two sets of the variables:
\bea
au+bv+cw&=&ax+by+cz=\xi,\\
 pu+qv+rw&=&px+qy+rz=\eta.
\eea
 In fact, now $\xi$ and $\eta$ are linearly dependent quantities,
\be\lb{3.16}
r\xi=c\eta.
\ee
In this situation, we may rewrite \eq{3.9} as
\be\lb{3.17}
u=x+\phi(\xi),\, v=y+\psi(\xi),\, w=z-\frac ac\phi(\xi)-\frac bc\psi(\xi),\, \xi=ax+by+cz=au+bv+cw,
\ee
where $\phi$ and $\psi$ are arbitrary polynomial functions of $\xi$, which is a direct 
3-dimensional extension of \eq{2.12} for which the inverse is obviously constructed as well:
\be
x=u-\phi(\xi),\quad y=v-\psi(\xi),\quad z=w+\frac ac\phi(\xi)+\frac bc\psi(\xi),\quad \xi=au+bv+cw,\quad c\neq0.
\ee
 Of course,
we now have $\deg(P)=\deg(P^{-1})$ and the equality in \eq{deg} never occurs in nontrivial situations where
$\min\{\deg(\phi),\deg(\psi)\}\geq2$.

We emphasize that the polynomial functions $\phi$ and $\psi$ in \eq{3.9} and \eq{3.17} are of arbitrary degrees in particular.

It may be of interest to explore a Monge--Amp\`{e}re equation type structure, as \eq{2.5} as we did in 2 dimensions for the Jacobian equation \eq{2.2}, for \eq{3.2}. For this purpose, we note that the first equation in \eq{3.3} implies that the vector $(f,g,h)$, being divergence free, has a vector potential,
$(A,B,C)$, satisfying
\be\lb{3.18}
(f,g,h)=\mbox{curl of } (A,B,C)=(C_y-B_z,A_z-C_x,B_x-A_y).
\ee
Hence \eq{3.2} becomes the following second-order nonlinear equation
\bea\lb{3.19}
&&\left|\begin{array}{cc}C_{xy}-B_{xz}&C_{yy}-B_{yz}\\A_{xz}-C_{xx}&A_{yz}-C_{xy}\end{array}\right|
+\left|\begin{array}{cc}A_{yz}-C_{xy}&A_{zz}-C_{xz}\\B_{xy}-A_{yy}&B_{xz}-A_{yz}\end{array}\right|
+\left|\begin{array}{cc}C_{xy}-B_{xz}&C_{yz}-B_{zz}\\B_{xx}-A_{xy}&B_{xz}-A_{yz}\end{array}\right|\nn\\
&&+\left|\begin{array}{ccc}C_{xy}-B_{xz}&C_{yy}-B_{yz}&C_{yz}-B_{zz}\\A_{xz}-C_{xx}&A_{yz}-C_{xy}&A_{zz}-C_{xz}\\B_{xx}-A_{xy}&B_{xy}-A_{yy}&B_{xz}-A_{yz} \end{array}\right|=0,
\eea
of a kind of the Hessian type, generalizing the Monge--Amp\`{e}re equation with an unknown vector field, instead of a scalar function. In view of \eq{3.18}, we can insert \eq{3.4} and
\eq{3.6} to obtain $(A,B,C)$ as a solution to \eq{3.19}. This procedure is known as taking the inverse of the curl operation: For ${\bf r}=(x,y,z)$, form the vector field
\be
{\bf F}(x,y,z)=(f,g,h)\times {\bf r}.
\ee
Then we have (cf. \cite{Sp})
\be
(A,B,C)(x,y,z)=\int_0^1 {\bf F}(tx,ty,tz)\dd t.
\ee

The vector potential $(A,B,C)$ here serves the role of the scalar potential $h$ in \eq{2.4} such that the system \eq{3.19} replaces \eq{2.5}. Thus, in light of Problem \ref{P2.1}, we state

\begin{problem}\lb{P3.1}
Can one obtain all polynomial solutions, in terms of the vector potential $(A,B,C)$, to the equation \eq{3.19} of the Monge--Amp\`{e}re type as stated?
\end{problem}

More generally, we may formulate the following problem directly related to the Jacobian problem in three dimensions.

\begin{problem}\lb{P3.2}
Find all polynomial solutions for the unknown $(f,g,h)$ to the first-order nonlinear equation \eq{3.2} of the Jacobian type.
\end{problem}

As another reduction of \eq{3.2}, we may set $h=H(f,g)$ where $H$ is a prescribed function of $f$ and $g$. Hence \eq{3.2} becomes
\be\lb{3.20}
f_x+g_y+H_f f_z+H_g g_z+J(f,g)(x,y)+H_g J(f,g)(x,z)+H_f J(g,f)(y,z)=0.
\ee

The under-determined equations \eq{3.19} and \eq{3.20} can be reduced further by imposing some appropriate constraints on the unknowns.

\section{General dimensions}\lb{sec4}
\setcounter{equation}{0}

In the general situation, with the notation
\be\lb{4.1}
(u_1,\dots,u_n)=P(x_1,\dots,x_n)=(x_1+f_1,\dots,x_n+f_n),\quad f_{i,j}=\frac{\pa f_i}{\pa x_j},\quad i,j=1,
\dots,n,
\ee
then it is clear that the Jacobian equation $J(P)=1$ or $\det(I +F)=1$ where $F=(f_{i,j})$ assumes the form
\be\lb{4.2}
E_1(F)+E_2(F)+\cdots+E_n(F)=0,
\ee
where $E_k(F)$ is the sum of all $k$ by $k$ principal minors of the matrix $F$, $k=1,\dots,n$, such that
$E_1(F)=\mbox{tr}(F)$ and $E_n(F)=\det(F)$ (cf. \cite{HJ}).

We now aim to obtain a family of solutions of \eq{4.2} of our interest that satisfy the Jacobian conjecture and extend what we found earlier in low dimensions.

For such a purpose and suggested by the study in Section 3, we use $(a_{ij})$ to denote an $(n-1)$ by $n$ matrix
in $\bfF$ and introduce the variables
\be\lb{4.3}
\xi_i=\sum_{j=1}^n a_{ij}x_j,\quad i=1,\dots,n-1.
\ee
Define
\be\lb{4.4}
u_j=x_j+f_j,\quad j=1,\dots,n,
\ee
where $f_1,\dots,f_{n-1}$ are arbitrary polynomials in $\xi_1,\dots,\xi_{n-1}$, respectively, but
\be\lb{4.5}
f_n=\sum_{j=1}^{n-1}b_j f_{j}(\xi_j),
\ee
where the coefficients $b_1,\dots,b_{n-1}\in\bfF$ are to be determined through the equation \eq{4.2} which
due to \eq{4.5}
is now slightly reduced into
\be\lb{x4.6}
E_1(F)+E_2(F)+\cdots+E_{n-1}(F)=0,
\ee
which is still rather complicated. For simplicity and in view of the study in Section 3, we
impose the following {\em full invariance} condition between the two sets of variables $x_1,\dots,x_n$ and $u_1,
\dots,u_n$:
\be\lb{4.6}
\xi_i=\sum_{j=1}^n a_{ij} x_j=\sum_{j=1}^n a_{ij}u_i,\quad i=1,\dots,n-1,
\ee
so that by virtue of \eq{4.4} we arrive at
\bea
\sum_{j=1}^n a_{ij} u_j&=&\sum_{j=1}^n a_{ij} x_j+\sum_{j=1}^{n-1} a_{ij}f_j +a_{in}f_n\nn\\
&=&\xi_i+\sum_{j=1}^{n-1}\left(a_{ij}+a_{in}b_j\right)f_j,\quad i=1,\dots,n-1,
\eea
which results in the solution
\be\lb{4.8}
b_j=-\frac{a_{ij}}{a_{in}},\quad a_{in}\neq0,\quad i,j=1,\dots,n-1.
\ee
This solution indicates that the quantities $\xi_1,\dots, \xi_{n-1}$ are linearly dependent:
\be
a_{in}\xi_{j}=a_{jn}\xi_{i},\quad i,j=1,\dots, n-1,
\ee
which extends \eq{3.16}. Since the functions $f_1,\dots,f_{n-1}$ are arbitrary, we may now set
\be\lb{4.11}
f_i=\phi_i(\xi),\quad i=1,\dots,n-1,\quad \xi=\sum_{j=1}^n a_j x_j.
\ee
Hence we obtain the polynomial map $P$ defined by
\be\lb{4.12}
u_1=x_1+\phi_1(\xi),\quad\dots,\quad u_{n-1}=x_{n-1}+\phi_{n-1}(\xi),\quad u_n=x_n-\sum_{i=1}^{n-1}\frac {a_i}{a_n}
\phi_i(\xi).
\ee
With \eq{4.12}, it is readily checked that the inverse of the polynomial map $P$ defined
in \eq{4.1} is
given by
\be\lb{4.13}
x_1=u_1-\phi_1(\xi),\quad\dots,\quad x_{n-1}=u_{n-1}-\phi_{n-1}(\xi),\quad x_n=u_n+\sum_{i=1}^{n-1}\frac{a_{i}}{a_{n}} \phi_i(\xi),
\ee
where $\phi_1,\dots,\phi_{n-1}$ are polynomial functions of the variable
$\xi=a_1u_1+\cdots+a_{n}u_{n}$. Of course we now have $\deg(P^{-1})=\deg(P)$.

Note that $\phi_1,\dots,\phi_{n-1}$, in nontrivial situations, consist of terms of degrees at least
2 of the variable $\xi$, which are arbitrary otherwise. Since
$DP(0)=I$, we automatically get $J(P)=1$. In particular, $(f_1,\dots,f_{n-1},f_n)$ so constructed
is a solution to the Jacobian equation \eq{4.2} such that the associated polynomial map $P$ given in
\eq{4.1} satisfies the Jacobian conjecture. 

As in the lower dimensional situations, the study here prompts us with the following problem.

\begin{problem}
Find all polynomial solutions for the unknown $(f_1,\dots,f_n)$ to the first-order Jacobian type nonlinear equation \eq{4.2}.
\end{problem}

Further reductions to \eq{4.2} may be carried out along the lines shown in Section 3 which may be used to construct other families of
polynomial solutions and are omitted here.

\section{Return to four dimensions}\lb{sec5}
\setcounter{equation}{0}

The  method of Section 3 may be used to obtain solutions involving partial invariance as well in higher dimensions but the computation become rather complicated. As an illustration, in this section, we
apply such a  method further to consider the four-dimensional situation.

For this purpose, consider the variables given in \eq{4.3}--\eq{4.5} with $n=4$. So we have
\be\lb{x5.1}
\xi_i=\sum_{j=1}^4 a_{ij}x_j,\quad u_i=x_i+f_i(\xi_i), \quad i=1,2,3;\quad u_4=x_4+f_4,\quad f_4=\sum_{i=1}^3 b_i f_i (\xi_i).
\ee
To solve \eq{x4.6}, we split it into three separate equations:
\be\lb{x5.2}
E_1=0,\quad E_2=0,\quad E_3=0.
\ee
From \eq{x5.1}, we have
\be
f_{i,j}=f_i'(\xi_i) a_{ij},\quad i=1,2,3, \quad j=1,2,3,4; \quad f_{4,j}=\sum_{i=1}^3 a_{ij}b_if'_i(\xi_i),\quad j=1,2,3,4,
\ee
which give us the expressions
\bea
E_1&=&(a_{11}+a_{14} b_1)f_1'+(a_{22}+a_{24}b_2)f'_2+(a_{33}+a_{34}b_3)f'_3,\lb{x5.4}\\
E_2&=&\left|\begin{array}{cc} a_{11}&a_{12}\\a_{21}&a_{22}\end{array}\right|f_1'f_2'+\left|\begin{array}{cc} a_{11}&a_{13}\\a_{31}&a_{33}\end{array}\right|f_1'f_3' +\left|\begin{array}{cc} a_{22}&a_{23}\\a_{32}&a_{33}\end{array}\right|f_2'f_3'+b_1 \left|\begin{array}{cc} a_{22}&a_{24}\\a_{12}&a_{14}\end{array}\right|f_1'f_2'+
b_1\left|\begin{array}{cc} a_{33}&a_{34}\\a_{13}&a_{14}\end{array}\right|f_1'f_3' \nn\\
&&+b_2\left|\begin{array}{cc} a_{11}&a_{14}\\a_{21}&a_{24}\end{array}\right|f_1'f_2' +b_2\left|\begin{array}{cc} a_{33}&a_{34}\\a_{23}&a_{24}\end{array}\right|f_2'f_3' 
+b_3\left|\begin{array}{cc} a_{11}&a_{14}\\a_{31}&a_{34}\end{array}\right|f_1'f_3'  +b_3\left|\begin{array}{cc} a_{22}&a_{24}\\a_{32}&a_{34}\end{array}\right|f_2'f_3' ,\lb{x5.5}\\
E_3&=&\left(\left|\begin{array}{ccc} a_{11}&a_{12}&a_{13}\\a_{21}&a_{22}&a_{23}\\a_{31}&a_{32}&a_{33}\end{array}\right|
+ b_1 \left|\begin{array}{ccc} a_{22}&a_{23}&a_{24}\\a_{32}&a_{33}&a_{34}\\a_{12}&a_{13}&a_{14}\end{array}\right|
+ b_2 \left|\begin{array}{ccc} a_{11}&a_{13}&a_{14}\\a_{31}&a_{33}&a_{34}\\a_{21}&a_{23}&a_{24}\end{array}\right|
+b_3\left|\begin{array}{ccc} a_{11}&a_{12}&a_{14}\\a_{21}&a_{22}&a_{24}\\a_{31}&a_{32}&a_{34}\end{array}\right|  \right)   f_1'f_2' f_3'.\lb{x5.6}
\eea
Since the functions $f_1,f_2,f_3$ are arbitrary, we are led by \eq{x5.2} and \eq{x5.4}--\eq{x5.6} to the equations
\bea
&&a_{11}+a_{14} b_1=0,\quad a_{22}+a_{24}b_2=0,\quad a_{33}+a_{34}b_3=0,\lb{x5.7}\\
&&\left|\begin{array}{cc} a_{11}&a_{12}\\a_{21}&a_{22}\end{array}\right|+b_1 \left|\begin{array}{cc} a_{22}&a_{24}\\a_{12}&a_{14}\end{array}\right|+
b_2\left|\begin{array}{cc} a_{11}&a_{14}\\a_{21}&a_{24}\end{array}\right|=0,\lb{x5.8}\\
&&\left|\begin{array}{cc} a_{11}&a_{13}\\a_{31}&a_{33}\end{array}\right|+b_1\left|\begin{array}{cc} a_{33}&a_{34}\\a_{13}&a_{14}\end{array}\right|
+b_3\left|\begin{array}{cc} a_{11}&a_{14}\\a_{31}&a_{34}\end{array}\right|=0,\lb{x5.9}\\
&&\left|\begin{array}{cc} a_{22}&a_{23}\\a_{32}&a_{33}\end{array}\right|
 +b_2\left|\begin{array}{cc} a_{33}&a_{34}\\a_{23}&a_{24}\end{array}\right|+
b_3\left|\begin{array}{cc} a_{22}&a_{24}\\a_{32}&a_{34}\end{array}\right|=0,\lb{x5.10} \\
&&
\left|\begin{array}{ccc} a_{11}&a_{12}&a_{13}\\a_{21}&a_{22}&a_{23}\\a_{31}&a_{32}&a_{33}\end{array}\right|
+ b_1 \left|\begin{array}{ccc} a_{22}&a_{23}&a_{24}\\a_{32}&a_{33}&a_{34}\\a_{12}&a_{13}&a_{14}\end{array}\right|
+ b_2 \left|\begin{array}{ccc} a_{11}&a_{13}&a_{14}\\a_{31}&a_{33}&a_{34}\\a_{21}&a_{23}&a_{24}\end{array}\right|+b_3\left|\begin{array}{ccc} a_{11}&a_{12}&a_{14}\\a_{21}&a_{22}&a_{24}\\a_{31}&a_{32}&a_{34}\end{array}\right|=0.\lb{x5.11}
\eea

To satisfy \eq{x5.7}, we assume $a_{14},a_{24},a_{34}\neq0$ to obtain
\be\lb{x5.12}
b_1=-\frac{a_{11}}{a_{14}},\quad b_2=-\frac{a_{22}}{a_{24}},\quad b_3=-\frac{a_{33}}{a_{34}}.
\ee
Inserting \eq{x5.12} into \eq{x5.8}--\eq{x5.10}, we have
\bea
&&(a_{12}a_{24}-a_{14}a_{22})(a_{11}a_{24}-a_{14}a_{21})=0,\lb{x5.13}\\
&&(a_{13}a_{34}-a_{14}a_{33})(a_{11}a_{34}-a_{14}a_{31})=0,\lb{x5.14}\\
&&(a_{23}a_{34}-a_{24}a_{33})(a_{22}a_{34}-a_{24}a_{32})=0.\lb{x5.15}
\eea
These equations extend \eq{3.7}. To fulfill \eq{x5.13}--\eq{x5.15} in a minimal manner, we impose, for example,
\bea
&&a_{12}a_{24}=a_{14}a_{22},\quad a_{13}a_{34}=a_{14}a_{33},\quad a_{23}a_{34}=a_{24}a_{33},\lb{x5.16}\\
&& a_{11}a_{24}\neq a_{14}a_{21},\quad a_{11}a_{34}\neq a_{14}a_{31},\quad a_{22}a_{34}\neq a_{24}a_{32}.\lb{x5.17}
\eea
It can be checked that, with \eq{x5.12} and \eq{x5.16}, the equation \eq{x5.11} is automatically satisfied. So we have thus obtained a solution to the differential equation of the problem.

With such a solution, we can construct the desired polynomial automorphism in four dimensions. To this goal, we again explore certain invariance structure as considered earlier. First, using
\eq{x5.1} and then \eq{x5.12} and \eq{x5.16}, we have
\bea
\sum_{j=1}^4 a_{1j}u_j&=&\sum_{j=1}^4 a_{1j}x_j+(a_{12}+a_{14}b_2)f_2+(a_{13}+a_{14}b_3)f_3\nn\\
&=&\xi_1+\left(a_{12}-\frac{a_{14} a_{22}}{a_{24}}\right)f_2+\left(a_{13}-\frac{a_{14} a_{33}}{a_{34}}\right)f_3=\xi_1.
\eea
That is, $\xi_1$ is an invariant. Similarly, we have
\bea
&&\sum_{j=1}^4 a_{2j}u_j=\xi_2+\left(a_{21}-\frac{a_{24} a_{11}}{a_{14}}\right)f_1,\lb{x5.19}\\
&&\sum_{j=1}^4 a_{3j}u_j=\xi_3+\left(a_{31}-\frac{a_{34} a_{11}}{a_{14}}\right)f_1+\left(a_{32}-\frac{a_{34} a_{22}}{a_{24}}\right)f_2,\lb{x5.20}
\eea
 indicating that $\xi_2,\xi_3$ are not invariants due to \eq{x5.17}. From these results, we immediately obtain the inverse map given by
\be\lb{x5.21}
x_i=u_i-f_i(\xi_i),\quad i=1,2,3,\quad x_4=u_4+\frac{a_{11}}{a_{14}}f_1(\xi_1)+\frac{a_{22}}{a_{24}}f_2(\xi_2)+\frac{a_{33}}{a_{34}}f_3(\xi_3),
\ee
where now
\bea
\xi_1&=&\sum_{j=1}^4 a_{1j}u_j,\lb{x5.22}\\
\xi_2&=&\sum_{j=1}^4 a_{2j}u_j-\left(a_{21}-\frac{a_{24} a_{11}}{a_{14}}\right)f_1(\xi_1),\lb{x5.23}\\
\xi_3&=&\sum_{j=1}^4 a_{3j}u_j-\left(a_{31}-\frac{a_{34} a_{11}}{a_{14}}\right)f_1(\xi_1)-\left(a_{32}-\frac{a_{34} a_{22}}{a_{24}}\right)f_2(\xi_2),\lb{x5.24}
\eea
iteratively in terms of the variables $u_i$'s.  From the construction, it is clear that
\bea
\deg(P)&=&\max\left\{\deg(f_1),\deg(f_2),\deg(f_3)\right\},\\
\deg(P^{-1})&=&\deg(f_1)\deg(f_2)\deg(f_3).\lb{x5.26}
\eea

Furthermore, the three nonminimal cases are also worth describing. One is when replacing \eq{x5.17} with
\be
 a_{11}a_{24}=a_{14}a_{21},\quad a_{11}a_{34}\neq a_{14}a_{31},\quad a_{22}a_{34}\neq a_{24}a_{32}.\lb{x5.27}
\ee
Using \eq{x5.27} in \eq{x5.19} and \eq{x5.20}, we see that the quantity $\xi_2$ becomes an additional invariant such that we update \eq{x5.23} into
\be\lb{x5.28}
\xi_2=\sum_{j=1}^4 a_{2j}u_j,
\ee
 but $\xi_3$ remains intact as a non-invariant. Hence the inverse map \eq{x5.21} is defined by \eq{x5.22}, \eq{x5.28}, and \eq{x5.24}, such that
\be\lb{x5.29}
\deg(P^{-1})=\deg(f_3)\max\{\deg(f_1),\deg(f_2)\}.
\ee
Another case is when replacing \eq{x5.27} with
\be
 a_{11}a_{24}=a_{14}a_{21},\quad a_{11}a_{34}= a_{14}a_{31},\quad a_{22}a_{34}\neq a_{24}a_{32},\lb{x5.30}
\ee
for example. Although the quantity $\xi_3$ is still a non-invariant, it is reduced into
\be
\xi_3=\sum_{j=1}^4 a_{3j}u_j-\left(a_{32}-\frac{a_{34} a_{22}}{a_{24}}\right)f_2(\xi_2),\lb{x5.31}
\ee
such that the inverse map defined by \eq{x5.21} with  \eq{x5.22}, \eq{x5.28}, and \eq{x5.31} satisfies
\be\lb{x5.32}
\deg(P^{-1})=\max\{\deg(f_1),\deg(f_2)\deg(f_3)\}.
\ee
The last case is when replacing \eq{x5.30} with all equalities,
\be
 a_{11}a_{24}=a_{14}a_{21},\quad a_{11}a_{34}= a_{14}a_{31},\quad a_{22}a_{34}= a_{24}a_{32},\lb{x5.33}
\ee
which finally renders the quantity $\xi_3$ an invariant as well,
\be\lb{x5.34}
\xi_3=\sum_{j=1}^4 a_{3j}u_j=\sum_{j=1}^4 a_{3j}x_j,
\ee
such that the inverse map \eq{x5.21} is defined by \eq{x5.22}, \eq{x5.28}, and \eq{x5.34}, with
\be\lb{x5.35}
\deg(P^{-1})=\max\{\deg(f_1),\deg(f_2),\deg(f_3)\}=\deg(P).
\ee
In fact, this last situation with a full set of invariant variables is contained in the general situation treated in Section 4.

The degree formulas \eq{x5.26}, \eq{x5.29}, \eq{x5.32}, and \eq{x5.35} are refined realizations of the general upper bound \eq{2.13} when $n=4$.

In summary of the study of this section, we state

\begin{theorem}
In the situation of four dimensions, the nonlinear first-order differential equation  \eq{4.2} with $n=4$ has a family of polynomial solutions depending
on three arbitrary polynomial functions, $f_1,f_2,f_3$, of the variables $\xi_1,\xi_2,\xi_3$ given as in \eq{x5.1}, respectively, where the coefficients $a_{ij}$'s satisfy
$a_{14},a_{24},a_{34}\neq0$ and the combined  conditions \eq{x5.12}--\eq{x5.15} for the existence of such solutions which give rise to
the inverse polynomial map explicitly through the expressions \eq{x5.21} so that one, two, or all three of the quantities $\xi_1,\xi_2,\xi_3$ appear as invariants
relating the two sets of the variables $x_i$'s and $u_i$'s.

\begin{enumerate}
\item[(i)] When one of the quantities $\xi_1,\xi_2,\xi_3$ is an invariant, the degree of the inverse map is the highest possible given by \eq{x5.26}.

\item[(ii)] When two of the quantities $\xi_1,\xi_2,\xi_3$ are invariants, say  $\xi_1$ and $\xi_2$, the degree of the inverse map is lowered to either \eq{x5.29} or
\eq{x5.32}, depending on the details of the residual non-invariant variable $\xi_3$ as described.

\item[(iii)] When all the quantities $\xi_1,\xi_2,\xi_3$ are invariants, the map and its inverse are described by the polynomials of the same structures as described, and in particular, the  map
and its inverse
have the same arbitrarily prescribed degrees as stated in \eq{x5.35}.
\end{enumerate}
\end{theorem}

The most transparent situation is when $f_1,f_2, f_3$ are of the same degree. With this assumption, the results \eq{x5.26}, \eq{x5.29}, \eq{x5.32}, and \eq{x5.35} read
\be
\deg(P^{-1})=(\deg(P))^3,\quad \deg(P^{-1})=(\deg(P))^2,\quad \deg(P^{-1})=\deg(P),
\ee
respectively, arranged in a descending order.

\section{The Jacobian problem with parametrization}\lb{sec6}
\setcounter{equation}{0}

The separation reductions from \eq{2.2} into \eq{2.3}, \eq{3.2} into \eq{3.3}, and \eq{x4.6}  (with $n=4$) into \eq{x5.2} prompt us to pursue a systematic and natural mechanism
to split the original Jacobian equation, $J(P)=1$, in order to render the problem effectively and schematically solvable, which is what we aim to do next.
In fact, it is clear that the Jacobian problem allows the following normalized reformulation
\be\lb{6.1}
(u_1,\dots,u_n)=P(x_1,\dots,x_n)=(\lm_1 x_1+f_1,\dots,\lm_n x_n+f_n),
\ee
generalizing \eq{4.1}, where $\lm_1,\dots,\lm_n$ are nonzero scalars, and $f_1,\dots,f_n$ are polynomials containing quadratic and higher power terms  in $x_1,\dots,x_n$. In this situation,
we have
\be\lb{6.2}
J(P)=\det((\lm_i \delta_{ij}+f_{i,j})), \quad f_{i,j}=\frac{\pa f_i}{\pa x_j},\quad i,j=1,\dots,n.
\ee
The condition on $f_1,\dots,f_n$ implies that $f_{i,j}$ all vanish at $x_1=0,\dots,x_n=0$, which leads to the equation 
\be\lb{6.3}
\det((\lm_i \delta_{ij}+f_{i,j}))=\lm_1\cdots\lm_n,
\ee
for the Jacobian problem. In order to explore this equation, we rewrite its left-hand side as
\be\lb{6.4}
\det((\lm_i \delta_{ij}+f_{i,j}))=C_0+\sum_{k=1}^{n-1}\sum_{\{i_1,\dots,i_k\}\subset\{1,\dots,n\},i_1<\cdots<i_k}C_{i_1\dots i_k}\lm_{i_1}\cdots\lm_{i_k}+\lm_1\cdots\lm_n,
\ee
where
\bea
C_0&=&\det((\lm_i \delta_{ij}+f_{i,j}))\big|_{\lm_1=\cdots=\lm_n=0}=E_n(F)=\det(F),\lb{6.4a}\\
C_{i_1\dots i_k}&=&\det(M_{i_1\cdots i_k})\big|_{\lm_{j_1}=\cdots=\lm_{j_{n-k}}=0,\{j_1,\dots,j_{n-k}\}=\{1,\dots,n\}\setminus\{i_1,\dots,i_k\}},\lb{6.4b}
\eea
where  $M_{i_1\cdots i_k}$ is the minor matrix of the Jacobian matrix given in \eq{6.2} with the $i_1,\dots,i_k$th rows and columns being deleted. 
In fact, \eq{6.4a} is obvious and \eq{6.4b} may be obtained from differentiating \eq{6.3} and \eq{6.4}
with respect to $\lm_{i_1},\dots,\lm_{i_k}$ and then setting $\lm_{j_1},\dots,\lm_{j_{n-k}}$ to zero,
where $\{j_1,\dots,j_{n-k}\}=\{1,\dots,n\}\setminus\{i_1,\dots,i_k\}$. As a consequence, we have
\be\lb{6.4c}
C_{i_1\dots i_k}=\det(F_{i_1\cdots i_k}),\quad \{i_1,\dots,i_k\}\subset\{1,\dots,n\},\quad i_1<\cdots<i_k,\quad k=1,\dots,n-1,
\ee
where  $F_{i_1\cdots i_k}$ is the minor matrix of the matrix $F=(f_{i,j})$ with the $i_1,\dots,i_k$th rows and columns being deleted.

In view of \eq{6.3} and \eq{6.4}, with \eq{6.4a}--\eq{6.4c}, we arrive at
\be\lb{6.4x}
\det(F)=0,\, \det(F_{i_1\dots i_k})=0,\, \{i_1,\dots,i_k\}\subset\{1,\dots,n\},\, i_1<\cdots<i_k,\, k=1,\dots,n-1.
\ee

It is clear that there are a total number of 
\be\lb{6.5}
\sum_{k=0}^{n-1}\left(\begin{array}{c} n\\k\end{array}\right)=2^n-1 
\ee
  equations in \eq{6.4c}, in contrast with a single equation in the original Jacobian problem. In other words, the reformulated
Jacobian problem enjoys a much tighter structure for the determination of its solution, and thus, may hopefully be resolved thoroughly.

We first consider the case in two dimensions so that \eq{2.1} assumes the form
\be\lb{6.6}
P(x,y)=(\alpha x+f(x,y),\beta y+g(x,y)),
\ee
where $\alpha,\beta$ are parameters.
Hence \eq{6.4x} gives rise to $2^2-1=3$ equations
\be\lb{6.7}
f_x=0,\quad g_y=0,\quad f_x g_y-f_y g_x=0.
\ee
The first two equations indicate that $f$ and $g$ depend on $y$ and $x$ only, respectively. So the last equation in \eq{6.7} is simply $f_y g_x=0$. That is $g=0$ or $f=0$.
Therefore we obtain two general solutions to \eq{6.7}:
\be\lb{6.8}
f=f(y),\quad g=0;\quad f=0,\quad g=g(x).
\ee
From these solutions and \eq{6.6}, we find under the condition $\alpha\neq0,\beta\neq0$  the corresponding polynomial maps and their inverses to be
\bea
&&u=\alpha x+f(y),\quad v=\beta y;\quad x=\frac 1\alpha\left(u-f\left(\frac v\beta\right)\right),\quad y=\frac{v}\beta,\\
&&u=\alpha x,\quad v=\beta y+g(x);\quad x=\frac u\alpha,\quad y=\frac1\beta\left(v-g\left(\frac u\alpha\right)\right),
\eea
respectively. We note that \eq{6.7} is a system of 3 equations which is under-determined for the 4 unknowns $f_x,f_y, g_x, g_y$ and its general solution \eq{6.8} contains 1 arbitrary function which accounts for 
the underlying extra degree of freedom, $4-3=1$.

We next consider the case in three dimensions so that \eq{3.1} becomes
\be\lb{6.11}
(u,v,w)=P(x,y,z)=(\alpha x+f(x,y,z),\beta y+g(x,y,z),\gamma z+h(x,y,z)).
\ee
Applying \eq{6.4x} to \eq{6.11}, we obtain $2^3-1=7$ equations
\bea
&&f_x=0,\quad g_y=0,\quad h_z=0,\lb{6.12}\\
&&f_y g_x=0,\quad f_z h_x=0,\quad g_zh_y=0,\lb{6.13}\\
&& f_zg_x h_y+f_y g_z h_x=0,\lb{6.14}
\eea
subsequently and iteratively.
The full family of general solutions to this system may be described with various combinations of the valid choices of $f,g,h$. In fact, \eq{6.12} indicates that $f=f(y,z),
g=g(x,z), h=h(x,y)$. That is, the three functions $f,g,h$ can only depend on up to two variables among $x,y,z$. Assume that $h$ depends on $x,y$. Then $h_x\neq0, h_y\neq0$.
Thus \eq{6.13} indicates that $f_z=g_z=0$. So $f,g$ can only depend on at most one variable nontrivially. We may pick $g$ to be such a polynomial so that $g$ depends on $x$. From
$g_x\neq0$, we obtain from the first equation in \eq{6.13} that $f_y=0$. That is, $f=0$. As a result, the equation \eq{6.14} is automatically satisfied. Hence we are lead to
$f=0$, $g=g(x)$, and $h=h(x,y)$, which give rise to the polynomial automorphism (the polynomial map and its polynomial inverse):
\bea
&&u=\alpha x,\quad v=\beta y+g(x),\quad w=\gamma z+h(x,y);\lb{6.15}\\
&&x=\frac u\alpha,\quad y=\frac1\beta\left(v-g\left(\frac u\alpha\right)\right),\quad z=\frac1\gamma\left(w-h\left(\frac u\alpha,\frac1\beta\left(v-g\left(\frac u\alpha\right)\right)\right)\right),\lb{6.16}
\eea
subject to $\alpha\neq0,\beta\neq0,\gamma\neq0$.
Now assume $h$ depends only on either $x$ or $y$, say $x$. Then $h_x\neq0$ and the second equation in \eq{6.13} gives us $f_z=0$. Using these in \eq{6.14}, we have
$f_y g_z=0$. If $f_y=0$, then $f=0$, and we arrive at the solution $f=0,g=g(x,z),h=h(x)$. If $f_y\neq0$, then $g_z=0$ and $g_x=0$ from the first equation in \eq{6.13}. Thus we arrive at the solution $f=f(y), g=0,h=h(x)$. 
Then assume that $h$ is independent of $x,y$. Hence $h=0$ and we have $f_y g_x=0$ from \eq{6.13} as the only remaining nontrivial equation. If $f_y=0$, then $f=f(z)$ and $g=g(x,z)$.
If $g_x=0$, then $f=f(y,z)$ and $g=g(z)$.
All of these solutions are of
the type \eq{6.15} after a variable relabeling.

Thus we conclude that one of the three polynomials $f, g, h$ must be trivial in general, and, modulo symmetry in the choice of variables, \eq{6.15} is the general polynomial
solution of the equations \eq{6.12}--\eq{6.14}.

In the case of four dimensions, we  consider the polynomial map $(u_1,\dots,u_4)=P(x_1,\dots,x_4)$ given by
\be\lb{6.17}
u_i=\lm_i x_i+f_i(x_1,\dots,x_4),\quad i=1,\dots,4.
\ee
With \eq{6.17}, the equations \eq{6.4x} become
\bea
&& f_{i,i}=0,\lb{6.18}\\
&&f_{i,j}f_{j,i}=0,\lb{6.19}\\
&&f_{i,j}f_{k,i}f_{j,k}+f_{j,i}f_{i,k}f_{k,j}=0,\lb{6.20}
\eea
where $i,j,k=1,\dots,4$, and 
\bea\lb{6.21}
&&f_{1,2}(f_{2,4}f_{3,1}f_{4,3}+f_{2,3}f_{3,4}f_{4,1})+f_{1,3}(f_{2,4}f_{3,2}f_{4,1}+f_{2,1}f_{3,4}f_{4,2})\nn\\
&&+f_{1,4}(f_{2,1}f_{3,2}f_{4,3}+f_{2,3}f_{3,1}f_{4,2})=0.
\eea
There are 4 equations in \eq{6.18}, 6 in \eq{6.19}, and 4 in \eq{6.20}, totaling $15=2^4-1$ equations in the coupled system \eq{6.18}--\eq{6.21}.  To search for the solution, we see from \eq{6.18} that each $f_i$ is
independent of $x_i$.  Furthermore, from \eq{6.19}, we have $f_{1,2}f_{2,1}=f_{1,3}f_{3,1}=f_{1,4}f_{4,1}=0$, which may be satisfied with setting $f_{1,i}=0$ for all $i$. That is,
$f_1=0$. This condition also renders \eq{6.21}. The rest of \eq{6.19} are
\be\lb{6.22}
f_{2,3}f_{3,2}=0,\quad f_{2,4}f_{4,2}=0,\quad f_{3,4}f_{4,3}=0,
\ee
which may be fulfilled with $f_2=f_2(x_1)$ and $f_3=f_3(x_1,x_2)$. With these results, it is readily examined that \eq{6.20} also follows automatically. This indicates that the equations
\eq{6.18}--\eq{6.21} are not all independent ones and that $f_4=f_4(x_1,x_2,x_3)$ is arbitrary. Consequently, we arrive at the following general solution of \eq{6.18}--\eq{6.21}:
\be\lb{6.23}
f_1=0,\quad f_2=f_2(x_1),\quad f_3=f_3(x_1,x_2),\quad f_4=f_4(x_1,x_2,x_3),
\ee
which gives rise to the polynomial automorphism
\bea
&&u_1=\lm_1 x_1,\quad u_2=\lm_2 x_2+f_2(x_1),\quad
 u_3=\lm_3 x_3+f_3(x_1,x_2),\quad u_4=\lm_4 x_4 +f_4(x_1,x_2,x_3);\\
&&x_1=\frac{u_1}{\lm_1},\quad x_2=\frac1{\lm_2}\left(u_2-f_2\left(\frac{u_1}{\lm_1}\right)\right),\quad
x_3=\frac1{\lm_3}\left(u_3-f_3\left(\frac{u_1}{\lm_1},\frac1{\lm_2}\left(u_2-f_2\left(\frac{u_1}{\lm_1}\right)\right)\right)\right),\nn\\
&&x_4=\frac1{\lm_4}\left(u_4-f_4\left(\frac{u_1}{\lm_1},\frac1{\lm_2}\left(u_2-f_2\left(\frac{u_1}{\lm_1}\right)\right), \frac1{\lm_3}\left(u_3-f_3\left(\frac{u_1}{\lm_1},\frac1{\lm_2}\left(u_2-f_2\left(\frac{u_1}{\lm_1}\right)\right)\right)\right)\right)\right),\nn\\
\lb{6.25}
\eea

Solutions associated with $f_2=0, f_3=0$, and $f_4=0$, respectively, are similarly constructed. 

We now show that the solutions of the above structure are all the solutions possible.

In fact, we first consider $f_{1,2}=0, f_{1,3}\neq0, f_{1,4}=0$. Thus $f_1=f_1(x_3)$ and $f_{3,1}=0$. The rest of
\eq{6.19} are again given by \eq{6.22}. We may choose $f_3=f_3(x_2)$ and $f_4=f_4(x_1,x_3)$ to fulfill these equations which lead to $f_{2,3}=f_{3,4}=0$. Inserting these results into
\eq{6.21}, we get $f_{1,3}f_{2,4}f_{3,2}f_{4,1}=0$ or $f_{2,4}f_{3,2}f_{4,1}=0$. We may choose $f_{2,4}=0$ for example. Then $f_2=f_2(x_1)$. On the other hand, from \eq{6.20}, we have
$f_{1,3}f_{k,1}f_{3,k}=0$, which gives us $f_{2,1}f_{3,2}=0$. If $f_{2,1}=0$, then $f_2=0$; if $f_{3,2}=0$, then $f_3=0$. Either case has been covered already as discussed.
We next consider $f_{1,2}=0, f_{1,3}\neq0, f_{1,4}\neq 0$. Thus $f_1=f_1(x_3,x_4)$ and $f_{3,1}=f_{4,1}=0$. So the rest of the equations in \eq{6.19} are still given by \eq{6.22}.
If $f_{2,3}\neq0$, then $f_{3,2}=0$. If $f_{4,2}=f_{4,3}=0$, then $f_4=0$ and we return to a familiar case. If $f_{4,2}=0$ but $f_{4,3}\neq0$, then $f_{3,4}=0$. Hence $f_3=0$ which is again
a familiar case.

The solution \eq{6.23} allows a direct generalization to any $n$-dimensional situation with
\be\lb{6.26}
f_1=0,\quad f_2=f_2(x_1),\quad f_3=f_3(x_1,x_2),\quad\dots\dots,\quad  f_n=f_n(x_1,\dots,x_{n-1}),
\ee
because the polynomial map given by
\be\lb{6.27}
u_1=\lm_1 x_1,\quad u_2=\lm_2 x_2+f_2(x_1),\quad \dots,\quad u_n=\lm_n x_n+f_n(x_1,\dots,x_{n-1}),
\ee
automatically satisfies the parametrized Jacobian equation, or
\be\lb{6.28}
\det\left(\begin{array}{ccccc}\lm_1&0&\cdots&\cdots&0\\f_{2,1}&\lm_2&\ddots&\ddots&\vdots\\\vdots&\ddots&\ddots&\ddots&\vdots\\\vdots&\vdots&\ddots&\ddots&0\\f_{n,1}&f_{n,2}&\cdots&f_{n,n-1}&\lm_n\end{array}\right)=\lm_1\lm_2\cdots\lm_n.
\ee
The inverse of \eq{6.27} is easily constructed as done in \eq{6.25} in four dimensions and thus omitted.

We now pursue a complete solution of the parametrized Jacobian equation \eq{6.3} following the ideas obtained in the above low-dimensional situations, in a more effective manner.

First,  the $n$-dimensional form of the equation \eq{6.18}  indicates that $f_i$ is independent of $x_i$ for $i=1,\dots,n$. Assume that $x_1$ (say) is a variable that enjoys a maximal presence in the problem so
that all $f_i$'s ($i\geq2$) depend on $x_1$. The $n$-dimensional form of the equation \eq{6.19} implies that if $f_i$ depends on $x_j$ then $f_j$ is independent of $x_i$, $i,j=1,\dots,n, i\neq j$.
Thus, with
\be
f_{2,1}\neq0,\quad\dots,\quad f_{n,1}\neq0,
\ee
we have $f_{1,2}=\cdots=f_{1,n}=0$. Hence $f_1$ and $f_2$ do not depend on $x_2$.  Assume that $x_2$ enjoys a maximal presence in the problem so that
all $f_i$'s ($i\geq3$) depend on $x_2$. Thus, with
\be
f_{3,2}\neq0,\quad\dots,\quad f_{n,2}\neq0,
\ee
we have $f_{2,3}=\cdots=f_{2,n}=0$. Carrying out this procedure to the $(n-1)$th step, we arrive at the conclusion that the polynomials $f_1,\dots, f_{n-1}$ do not depend
on $x_{n-1}$. So we are left with assuming that $f_n$ depends on $x_{n-1}$ to achieve a maximal presence of the variable $x_{n-1}$ in the set of the polynomials. This gives us
$f_{n,n-1}\neq0$. Hence $f_{n-1,n}=0$. Consequently, this procedure leads us to \eq{6.28} and the solution \eq{6.27}, constructively and systematically.

Alternatively, if we assume the polynomial function $f_1$ enjoys the maximal presence of the variables, then $f_1=f_1(x_2,\dots,x_n)$. This assumption leads to $f_{1,2}\neq0,\dots,
f_{1,n}\neq0$ and $f_{2,1}=\cdots=f_{n,1}=0$. Subsequently, assume $f_2$ enjoys the maximal presence of the rest of the variables. Since $f_{2,1}=0$, we have
$f_2=f_2(x_3,\dots,x_n)$. Following this process, we eventually have $f_{n-1}=f_{n-1}(x_n)$ and $f_n=0$. That is,
\be\lb{x6.31}
f_1=f_1(x_2,x_3,\dots,x_n),\quad f_2=f_2(x_3,\dots,x_n),\quad\dots,\quad f_{n-1}=f_{n-1}(x_n),\quad f_n=0.
\ee
This solution differs from \eq{6.26} only by a relabeling of the variables, of course.

This result in fact also includes all the degenerate cases. For example, assume only $f_{1,2}\neq0,\dots,f_{1,n-1}\neq0$ but $f_{1,n}=0$. That is, $x_n$ does not make an
explicit presence in $f_1$. Then we have 
\be\lb{x6.32}
f_{2,1}=0,\quad\dots,\quad f_{n-1,1}=0,
\ee
and we have either $f_{n,1}=0$ or $f_{n,1}\neq0$. In the former case, the equation \eq{6.3} becomes
\be\lb{x6.33}
\det\left(\begin{array}{cccc}\lm_2&f_{2,3}&\cdots&f_{2,n}\\f_{3,2}&\cdots&\cdots&f_{3,n}\\\vdots&\ddots&\ddots&\vdots\\ f_{n-1,2}&\cdots&\ddots&f_{n-1,n}\\f_{n,2}&f_{n,3}&\cdots&\lm_n\end{array}\right)=\lm_2\cdots\lm_n.
\ee
Using induction, we see that, with a relabeling of the variables if necessary, the general solution of \eq{x6.33}, including all degenerate cases, reads
\be
f_2=f_2(x_3,x_4,\dots,x_n),\quad f_3=f_2(x_4,\dots,x_n),\quad\dots,\quad f_{n-1}=f_{n-1}(x_n),\quad f_n=0,
\ee
following \eq{x6.31}. In view of \eq{x6.32}, $f_{n,1}=0$, and \eq{x6.33}, we obtain $f_1=f_1(x_2,\dots,x_{n-1})$, which is consistent with \eq{x6.31}. In the latter case, that is,
$f_{n,1}\neq0$, then the $n$-dimensional version of \eq{6.20} with $i=n,j=1$ gives us
\be
f_{k,n}f_{1,k}=0,\quad k=1,\dots,n.
\ee
However, the assumption $f_{1,k}\neq0$ for $k=2,\dots,n-1$ implies $f_{k,n}=0$ for $k=2,\dots,n-1$. Therefore $f_1,\dots,f_{n-1},f_n$ are independent of $x_n$ and the problem is
actually reduced into an $(n-1)$-dimensional problem. By induction, we conclude that the solution structure \eq{x6.31} again appears with the absence of an explicit dependence on $x_n$.

It should be noted that the general solution \eq{6.27} is not of the homogeneous type studied earlier. Besides, with fixed $\lm_1,\dots,\lm_n$, say $\lm_1=\cdots=\lm_n=1$, the result \eq{6.27} gives rise to a new family
of solutions to the classical Jacobian equation $J(P)=1$, say, as well.

In summary, we state the results of the study of this section as follows.

\begin{theorem}
In any $n\geq2$ dimensions, the general solution to the Jacobian equation
$J(P)=\lm_1\cdots\lm_n$ associated with the polynomial map \eq{6.1} parametrized by the parameters $\lm_1,\dots,\lm_n$
is given by the expression \eq{6.26}, up to a relabeling of the variables. For specific values of these parameters, 
the obtained solution renders the associated map invertible
with a polynomial inverse if and only if
the Jacobian of the map is nonvanishing which is characterized by $\lm_1\neq0,\dots,\lm_n\neq0$. In this situation,  the inverse map is iteratively given by
\be\lb{6.31}
x_1=\frac{u_1}{\lm_1},\quad x_2=\frac1{\lm_2}\left(u_2-f_2(x_1)\right),\quad \dots,\quad x_n=\frac1{\lm_n}\left(u_n-f_n(x_1,\dots,x_{n-1})\right),
\ee
schematically and explicitly.
In particular, this construction provides a new family of solutions to the
non-parametrized (classical) Jacobian problem, with normalized parameters, $\lm_1=\cdots=\lm_n=1$,
which is of non-homogeneous type in general.
\end{theorem}

From \eq{6.31}, it is clear that the degree of the inverse map $P^{-1}$ satisfies the upper bound
\be
\deg(P^{-1})\leq \deg(f_1)\cdots\deg(f_{n-1})\leq \deg(P)^{n-1},
\ee
which reconfirms the bound \eq{2.13} again.
\medskip

Of independent interest is the situation when the multiple parameters $\lm_1,\dots,\lm_n$ in \eq{6.1} are combined into a single one, $\lm_1=\cdots=\lm_n=\lm$. Then the Jacobian
equation \eq{6.3} becomes
\be\lb{a6.26}
\det(\lm I+F)=\lm^n,\quad F=(f_{i,j}),
\ee
such that the left-hand side of this equation is the characteristic polynomial of the matrix $-F$. As a consequence, we are led from \eq{a6.26} to \cite{HJ}:
\be\lb{a6.27}
E_1(F)\lm^{n-1}+E_2(F)\lm^{n-2}+\cdots+E_{n-1}(F)\lm+E_n(F)=0,
\ee
or equivalently,
\be\lb{a6.28}
E_i(F)=0,\quad i=1,\dots,n.
\ee
That is, the presence of the parameter $\lm$ renders a natural splitting of the equation \eq{4.2} into $n$ subequations, listed in \eq{a6.28}.

When $n=2,3, 4$, \eq{a6.28} becomes \eq{2.3}, \eq{3.3}, \eq{x5.2} (along with $E_4=0$), respectively, whose solutions  have been constructed in Sections 2, 3, 5. In those situations,
the equations are {\em imposed} in order to reduce the original Jacobian equation. In the present situation, these equations are natural consequences of the presence of the free parameter
$\lm$. However, unlike in the multiply parametrized situation, the singly parametrized system of the equations, \eq{a6.28}, defies a complete solution even in the bottom-dimensional 
situation, $n=2$, which explains the associated difficulties in the original problem transparently. 

At this juncture, it will be interesting to work out a concrete example as a direct illustration of our study. For simplicity, we consider the two-dimensional case \eq{2.1} with
\be
f(x,y)=\sum_{i+j=k,k=2,3} a_{ij}x^i y^j,\quad g(x,y)=\sum_{i+j=k,k=2,3} b_{ij}x^i y^j,\lb{6.35}\\
\ee
Substituting \eq{6.35} into the full two-dimensional Jacobian equation \eq{2.2} and solving for coefficients, we get four families of solutions. The two simple ones are given by
\be
f(x,y)=a y^2+by^3,\quad g(x,y)=0;\quad f(x,y)=0,\quad g(x,y)= ax^2+bx^3,
\ee
where $a,b$ are arbitrary scalars. These solutions are special cases contained already in \eq{6.8}. The other two families of solutions are given by
\bea
&&f(x,y)=a\xi^2-\frac{b^2}{9c}\xi^3,\quad g(x,y)=-\frac{ab}{3c}\xi^2+\frac{b^3}{27c^2}\xi^3,\quad \xi=x+\frac{3c}b y,\quad b\neq0,c\neq0,\lb{6.37}\\
&&f(x,y)=\frac{\eta^2}a,\quad g(x,y)=\frac{b\eta^2}{a^2},\quad \eta=bx-ay,\quad a\neq0,\lb{6.38}
\eea
in terms of arbitrary scalars, $a,b,c$, again. In these families of solutions, the variables $\xi$ and $\eta$ are invariant variables, through $u=x+f,v=y+g$. That is, we have
\be
\xi= x+\frac{3c}b y=u+\frac{3c}b v,\quad \eta=bx-ay=bu-av,
\ee
respectively, enabling immediate constructions of inverse polynomial maps, as a consequence. These solutions are special cases of what obtained in Section \ref{sec2}.  It is also clear that all the solutions here are ``divergence free", satisfying the
equation $f_x+g_y=0$, and thus the homogeneous Monge--Amp\`{e}re equation as well.

Besides, note that, in three dimensions, the system \eq{a6.28} splits the Hessian equation \eq{3.19} into a system of two separate equations:
\bea
&&\left|\begin{array}{cc}C_{xy}-B_{xz}&C_{yy}-B_{yz}\\A_{xz}-C_{xx}&A_{yz}-C_{xy}\end{array}\right|
+\left|\begin{array}{cc}A_{yz}-C_{xy}&A_{zz}-C_{xz}\\B_{xy}-A_{yy}&B_{xz}-A_{yz}\end{array}\right|
+\left|\begin{array}{cc}C_{xy}-B_{xz}&C_{yz}-B_{zz}\\B_{xx}-A_{xy}&B_{xz}-A_{yz}\end{array}\right|=0,\lb{6.50}\\
&&\left|\begin{array}{ccc}C_{xy}-B_{xz}&C_{yy}-B_{yz}&C_{yz}-B_{zz}\\A_{xz}-C_{xx}&A_{yz}-C_{xy}&A_{zz}-C_{xz}\\B_{xx}-A_{xy}&B_{xy}-A_{yy}&B_{xz}-A_{yz} \end{array}\right|=0,\lb{6.51}
\eea
and a polynomial integrability problem for these coupled equations may be formulated in terms of the vector unknown $(A,B,C)$, similar to
Problems \ref{P3.1} and \ref{P3.2}.

\section{Application to cryptography}\lb{sec7}
\setcounter{equation}{0}

Due to their analytic simplicity,  invertible polynomial maps are of obvious interest and importance in
applications. For example, 
we consider an application of the construction here to cryptography \cite{Hoff,Kn,NC} in which a polynomial automorphism can be used to serve as a cryptographic key to encrypt and decrypt messages. 
To see how, we note that although our construction of inverse maps in various dimensions is suggested by the Jacobian problem and the associated nonlinear partial differential equations, the 
 detailed fine structures of such maps are quite general. In other words, the final formulas obtained of the inverse maps by-pass the original technical assumptions and formulation of the problem.
For example, the construction works equally well when the underlying field $\bfF$ is of a nonzero characteristic or replaced by a commutative ring with a unity, say $\cal R$.
As a concrete example,  assume that a plaintext is arranged in the form of a sequence of vectors in ${\cal R}^n$, denoted by $x=(x_1,\dots,x_n)$. As a consequence of
the construction comprised of \eq{4.11} and \eq{4.12}, we may formulate the (nonlinear and polynomial) encryption map $u=(u_1,\dots,u_n)\equiv E(x)\in{\cal R}^n$ with
\bea
&&u_1=x_1+\phi_1(\xi),\quad\dots,\quad u_{n-1}=x_{n-1}+\phi_{n-1}(\xi),\quad u_n=x_n-\sum_{i=1}^{n-1} a_i\phi_i(\xi),\lb{a6.12}\\
&&\xi=\sum_{i=1}^{n-1} a_i x_i+ x_n,\lb{a6.13}
\eea
where $a_1,\dots,a_{n-1}\in {\cal R}$ are arbitrary ring elements and $\phi_1(\xi),\dots,\phi_{n-1}(\xi)$ are arbitrary polynomial functions of $\xi$. In this way, we come up with a
ciphertext. Note that, we choose $a_n=1$ in \eq{4.11} in order to avoid division by it in \eq{4.12}. Therefore, by \eq{4.13}, we get the decryption map
$x=D(u)$ with
\bea
&&x_1=u_1-\phi_1(\xi),\quad\dots,\quad x_{n-1}=u_{n-1}-\phi_{n-1}(\xi),\quad x_n=u_n+\sum_{i=1}^{n-1}{a_{i}} \phi_i(\xi),\lb{a6.14}\\
&&\xi=\sum_{i=1}^{n-1} a_i u_i+ u_n.\lb{a6.15}
\eea
 The invariance of the quantity $\xi$ relating the plaintext and ciphertext variable vectors, $x$ and $u$, respectively, is self-evident. The pair of maps, 
\be
E,D: {\cal R}^n\to{\cal R}^n,
\ee
collectively given in \eq{a6.12}--\eq{a6.15}, 
is seen to provide a  public-key cryptography, in which the key $\cal K$ could simply be taken to be a set of the vectors $(a_1,\dots,a_{n-1})$ in ${\cal R}^{n-1}$ consisting of the
coefficients of the map \eq{a6.12}--\eq{a6.13}, and more generally and complicatedly, the set of the same vectors in ${\cal R}^{n-1}$ supplemented with  the $n-1$ {\em arbitrary} polynomials $\phi_1,\dots,
\phi_{n-1}$ nonlinear in the variable $\xi$ as stated.

Moreover, to achieve further enhanced or elevated security, we may also use the polynomial automorphisms with partial invariance properties as constructed in Sections \ref{sec3} and \ref{sec5},
or those without any invariance structure presented in Section \ref{sec6}, as well as the more generalized ones to be given in the next section.

\section{The parametrized Jacobian problem in a general setting}\lb{s8}
\setcounter{equation}{0}

For broader applicability, we reformulate the problem in a general setting such that the multiparametrization is clearly represented and exhibited. In this context, the map \eq{1.4} may be rewritten as
\be\lb{8.1}
u=P(x)=Lx+N(x),\quad x,u\in\bfF^n,
\ee
where $L$ is an $n\times n$ matrix over $\bfF$ and $N:\bfF^n\to \bfF^n$ denotes the nonlinear part of the map. Then the Jacobian condition becomes
\be
J(P)=J(P)(0)=\det(L)\neq0.
\ee
With the new variable $y=Lx=(y_1,\dots,y_n)$, the map assumes its normalized form
\be
u=y +My,\quad M(x)=N(L^{-1} y),
\ee
such that the multiparametrized problem becomes
\bea
u&=&\Lambda y+M(y),\quad \Lambda=\mbox{diag}\{\lm_1,\dots,\lm_n\}\nn\\
&=& \Lambda Lx+N(x).
\eea
In this way, the Jacobian equation is recast into
\be\lb{8.5}
J(P)=\det(\Lambda L)=\lm_1\cdots\lm_n\det(L).
\ee

Our study has established the following theorem regarding the generalized multiparametrized Jacobian equation \eq{8.5}.

\begin{theorem}
Consider an $n$-dimensional polynomial map of the multiparametrized form $u=P(x)=\Lambda Lx+N(x)$ that extends the map \eq{8.1}. If this map has a polynomial inverse, then its Jacobian satisfies the equation \eq{8.5} with
\be\lb{8.6}
\lm_1\neq0,\quad \dots,\quad\lm_n\neq0,\quad \det(L)\neq0.
\ee
Conversely, subject to the condition \eq{8.6} and with the notation $L=(L_{ij})$, the Jacobian equation \eq{8.5} can be integrated over the space of polynomials whose general solution is of the normalized form
\be
N(x)=\left(0,f_2\left(\sum_{j=1}^n L_{1j}x_j\right),\dots,f_n\left(\sum_{j=1}^n L_{1j}x_j,\dots, \sum_{j=1}^n L_{n-1,j}x_j\right)\right),
\ee
where $f_2,\dots,f_n$ are arbitrary polynomials consisting of terms of nonlinear powers of the variables, such that the associated inverse map is given by the formulas
\bea
x=L^{-1}y,\quad&& y_1=\frac1{\lm_1}u_1,\quad y_2=\frac1{\lm_2}\left(u_2-f_2\left(\sum_{j=1}^n L_{1j}x_j\right)\right),\nn\\
&& \dots, \quad y_n=\frac1{\lm_n}\left(u_n-f_n\left(\sum_{j=1}^n L_{1j}x_j,\dots, \sum_{j=1}^n L_{n-1,j}x_j\right)\right),
\eea
as expressed in terms of an iterative procedure.
\end{theorem}

\section{Summary and comments}\lb{sec8}
\setcounter{equation}{0}

In this work, the Jacobian conjecture is directly studied by way of the partial differential equations
it prompts. These equations are first order, nonlinear, and expressed in terms of the sum of all
principal minors of the Jacobian matrix of the nonlinear part of the polynomial map. In two dimensions, the
equation may be reduced into a homogeneous Monge--Amp\`{e}re equation whose solutions uncover an invariance
structure between the variables which may be exploited effectively beyond two dimensions. As a consequence, in any general $n$ dimension ($n\geq2$),  solutions
depending on $n-1$ arbitrary polynomial functions are obtained, which give rise to polynomial maps
of arbitrarily high degrees, satisfying the
Jacobian conjecture. 
The maps may or may not be homogeneous, depending on the choice of
these arbitrary polynomial functions. Interestingly and practically, for these maps, the inverse maps are of similar structures and
can be constructed immediately using the established invariance property,
either partial or full, relating the two sets of variables in consideration. 
More importantly, it is demonstrated that the Jacobian problem enjoys a multiply parametrized reformulation which splits the Jacobian equation into a system of coupled equations that may be
integrated effectively and explicitly such that the nonzero Jacobian condition ensures that the polynomial map is invertible and
possesses a polynomial inverse map. This inverse map may be constructed schematically in which no use is made of any invariance property. This construction
 settles the multiply parametrized Jacobian problem in its general setting.

\medskip

{\bf Data availability statement:} The data that support the findings of this study are available within the article.

{

}

\end{document}